\documentclass[12pt]{article}
\usepackage{amsmath,amssymb,amsthm}
\usepackage{boxedminipage}
\usepackage{graphicx}
\usepackage[margin=1in]{geometry}
\usepackage{mathpazo}
\usepackage{fancyhdr}
\usepackage{enumitem}
\usepackage{longtable}

\usepackage{tikz}
\usetikzlibrary{arrows}

\title{\Large{Counterexamples to a conjecture by Alaoglu and Erd\H{o}s about superabundant numbers}}
\author{Burdette, Tibor \& Stewart, Ian M.}

\newtheorem{conj}{Conjecture}

\begin{document}

%\collaborators{ME'N'EEN}

%%% PROBLEMS AND SOLUTIONS -------------------------------------------------------------------
% Include all problems statements
% For solutions, use the \begin{proof} environment
% Note that \begin{proof}[Solution] will replace the label "Proof" with "Solution", etc.
%
\maketitle
Superabundant (SA) numbers $\mathcal{S}$ are the set of integers $n$ for which \[\frac{\sigma(n)}{n} > \frac{\sigma(m)}{m}\text{ for }1 \le m < n,\] where $\sigma(n)$ is the sum of the divisors of $n$.
 
In 1944, Alaoglu and Erd\H{o}s\cite{alaoglu_erdos} conjectured that for every superabundant number $n$, there exist primes $p,q$ such that the product $np$ and the quotient $n/q$ are superabundant.

Another way to state this conjecture is that the infinite lattice of superabundant numbers (where $n_1$ is linked to $n_2$ if $pn_1=n_2$ for some prime $p$) has one source (the number 1) and no sinks.

Assuming the truth of this conjecture, along with the property that the prime decomposition of any $n\in \mathcal{S}$ has nonincreasing exponents\cite{alaoglu_erdos}, an algorithm can be constructed that inductively generates all SA numbers from 1 up to any arbitrary bound. Cross-referencing the results of this algorithm with T. D. Noe's\cite{noe} list of the first million SA numbers reveals the existence of numbers $n$ in $\mathcal{S}$ not produced by the algorithm - i.e., there does not exist $q$ s.t. $n/q\in\mathcal{S}$. The smallest of these, using the notation defined in the appendix, is \[\{738,27,8,5,4,3,0,0,2,0,0,0,0,0,1\}\approx 6.04\times10^{2448},\] index 19861 in Noe's list.

Searching through the Noe dataset reveals 18 counterexamples to the $n/q$ half of Alaoglu and Erd\H{o}s's conjecture (sources) and 88 counterexamples to the $np$ half of Alaoglu and Erd\H{o}s's conjecture (sinks) less than $10^{100,000}$. These are listed in the Appendix.

We propose the following weaker version of Alaoglu and Erd\H{o}s's conjecture:
        \begin{conj}
        	The infinite lattice of superabundant numbers is connected. That is, given a superabundant number n:
        \begin{enumerate}
        \item it is possible to find either $p$ such that $np$ is superabundant or $q$ such that $n/q$ is superabundant, and
        \item it is always possible to find a chain of primes $p_1 p_2 p_3 ... p_m$ such that by successive multiplication OR division of these primes it is possible to reach 1.
        \end{enumerate}
        \end{conj}
        If this conjecture is proven to be true, then it would be possible to generate the superabundant numbers using a modified version of the previous algorithm.

\newpage
\section*{Appendix: List of all counterexamples less than $10^{100,000}$}
The numbers in the table below are much too large to express in standard form, so we represent them in a "super-compactification notation" (SCN) defined as follows:

\textit{If $p$ is in position $k$ when $n$ is represented in SCN, $p\#$ divides $n$ exactly $k$ times.}\\
A few examples are listed below:
\begin{itemize}
    \item $\{1\}=2$
    \item $\{0,0,1\}=2^3$
    \item $\{3\}=2\cdot3\cdot5=30$
    \item $\{4,0,1\}=2^3\cdot3\cdot5\cdot7=840$
\end{itemize}
Indices are from Noe's dataset.
\begin{longtable}{r|c|r|r|l}\hline
\textbf{Index}&\textbf{Type}&\textbf{Group}&$\boldsymbol{\log_{10} n}$&\textbf{SCN Representation}\\\hline\hline\endhead
2687&$np$&184&389.55&\{152,12,5,3,0,0,2,0,0,0,0,0,1\}\\
5780&$np$&337&849.25&\{295,16,7,4,3,0,0,2,0,0,0,0,0,1\}\\
5804&$np$&338&852.72&\{296,16,7,4,3,0,0,0,2,0,0,0,1\}\\
6180&$np$&354&906.34&\{312,17,6,4,0,3,0,2,0,0,0,0,1\}\\
8528&$np$&446&1209.21&\{400,19,7,4,0,3,0,0,2,0,0,0,1\}\\
9721&$np$&487&1344.58&\{438,20,8,0,4,3,0,2,0,0,0,0,0,1\}\\
17859&$np$&733&2206.25&\{674,24,9,5,4,3,0,0,2,0,0,0,0,0,0,0,1\}\\
19861&$n/q$&797&2448.78&\{738,27,8,5,4,3,0,0,2,0,0,0,0,0,1\}\\
20177&$np$&808&2482.5&\{747,26,9,5,4,3,0,0,2,0,0,0,0,0,0,0,1\}\\
22010&$np$&864&2692.77&\{802,27,9,6,4,3,0,0,2,0,0,0,0,0,0,1\}\\
22843&$np$&892&2793.19&\{828,28,9,5,4,3,0,0,0,2,0,0,0,0,0,0,1\}\\
24690&$n/q$&965&3082.71&\{903,30,9,5,4,3,0,2,0,0,0,0,0,0,1\}\\
29464&$np$&1144&3763.88&\{1077,30,10,6,4,3,0,0,2,0,0,0,0,0,0,0,1\}\\
29793&$np$&1154&3803.28&\{1087,31,9,5,4,3,0,0,0,2,0,0,0,0,0,0,1\}\\
34881&$np$&1297&4363.81&\{1227,32,10,6,4,3,0,0,0,2,0,0,0,0,0,0,1\}\\
34894&$np$&1297&4364.82&\{1227,32,11,6,4,3,0,0,2,0,0,0,0,0,0,0,1\}\\
40980&$np$&1466&5035.05&\{1392,34,11,6,4,3,0,0,0,2,0,0,0,0,0,0,0,1\}\\
41177&$np$&1471&5058.35&\{1398,34,10,6,4,0,3,0,0,2,0,0,0,0,0,0,1\}\\
42002&$np$&1495&5153.74&\{1421,35,10,6,4,3,0,0,0,2,0,0,0,0,0,0,0,1\}\\
43158&$np$&1533&5306.19&\{1458,35,11,6,4,3,0,0,0,2,0,0,0,0,0,0,0,1\}\\
43306&$np$&1537&5325.53&\{1463,35,10,6,4,0,3,0,0,2,0,0,0,0,0,0,1\}\\
46397&$np$&1647&5768.66&\{1570,36,11,6,4,0,3,0,0,2,0,0,0,0,0,0,0,1\}\\
47551&$n/q$&1692&5958.98&\{1616,37,10,6,4,0,3,0,0,2,0,0,0,0,0,0,1\}\\
50451&$np$&1797&6396.97&\{1721,37,12,6,4,3,0,0,0,2,0,0,0,0,0,1\}\\
60339&$np$&2139&7821.5&\{2057,41,12,6,4,3,0,0,0,2,0,0,0,0,0,0,0,1\}\\
64109&$np$&2260&8331.59&\{2176,42,11,6,0,4,3,0,0,2,0,0,0,0,0,0,0,1\}\\
66630&$np$&2328&8621.31&\{2243,43,11,6,0,4,3,0,0,2,0,0,0,0,0,0,0,1\}\\
66653&$np$&2329&8623.26&\{2243,43,12,7,4,0,3,0,0,2,0,0,0,0,0,0,0,1\}\\
69557&$np$&2397&8913.92&\{2310,44,12,7,4,0,3,0,0,2,0,0,0,0,0,0,0,1\}\\
71020&$np$&2433&9069.23&\{2346,45,11,6,4,0,3,0,0,0,2,0,0,0,0,0,0,1\}\\
112028&$np$&3401&13303.7&\{3302,51,14,7,5,4,3,0,0,0,2,0,0,0,0,0,0,1\}\\
120134&$np$&3605&14219.83&\{3505,53,13,7,0,4,3,0,0,0,0,2,0,0,0,0,0,1\}\\
127141&$np$&3780&15014.36&\{3680,53,15,7,0,4,3,0,0,0,2,0,0,0,0,0,1\}\\
135611&$np$&3967&15855.57&\{3864,55,14,7,0,4,0,3,0,0,2,0,0,0,0,0,0,1\}\\
143374&$np$&4112&16515.85&\{4008,56,14,7,0,4,3,0,0,0,0,2,0,0,0,0,0,1\}\\
143413&$np$&4114&16518.62&\{4008,56,15,8,0,4,3,0,0,0,2,0,0,0,0,0,0,0,1\}\\
143910&$n/q$&4120&16558.54&\{4017,57,14,7,0,4,3,0,0,2,0,0,0,0,0,0,0,1\}\\
148924&$np$&4219&16997.94&\{4112,57,15,8,0,4,3,0,0,0,2,0,0,0,0,0,0,0,1\}\\
150197&$n/q$&4239&17099.36&\{4134,58,15,7,0,4,3,0,0,2,0,0,0,0,0,0,0,1\}\\
155995&$n/q$&4350&17607.86&\{4244,59,15,7,0,4,3,0,0,2,0,0,0,0,0,0,0,1\}\\
157248&$np$&4378&17727.7&\{4270,58,15,8,0,4,0,3,0,0,2,0,0,0,0,0,0,1\}\\
158989&$n/q$&4415&17900.67&\{4307,59,15,8,5,4,3,0,0,2,0,0,0,0,0,0,0,1\}\\
161559&$np$&4473&18160.42&\{4363,59,15,8,5,4,0,3,0,0,2,0,0,0,0,0,0,1\}\\
167169&$n/q$&4603&18767.8&\{4494,61,14,7,5,4,3,0,0,0,2,0,0,0,0,0,0,1\}\\
187973&$np$&5151&21301.74&\{5036,62,15,8,5,4,0,3,0,0,2,0,0,0,0,0,0,0,0,1\}\\
188749&$n/q$&5167&21391.58&\{5055,63,15,8,5,4,3,0,0,2,0,0,0,0,0,0,0,1\}\\
193772&$np$&5296&21987.34&\{5182,62,15,8,5,4,0,3,0,0,0,2,0,0,0,0,0,1\}\\
218788&$np$&5943&25036.66&\{5824,66,16,8,5,4,3,0,0,0,2,0,0,0,0,0,0,0,0,1\}\\
219403&$np$&5956&25104.47&\{5838,66,17,8,5,4,3,0,0,0,2,0,0,0,0,0,0,1\}\\
238625&$np$&6400&27212.15&\{6278,68,16,9,5,4,3,0,0,0,2,0,0,0,0,0,0,0,0,1\}\\
241557&$np$&6467&27538.57&\{6346,68,16,9,5,4,0,3,0,0,2,0,0,0,0,0,0,1\}\\
242962&$np$&6502&27698.84&\{6379,69,16,8,5,4,0,3,0,0,2,0,0,0,0,0,0,0,0,1\}\\
246157&$np$&6572&28040.85&\{6450,69,16,9,5,4,0,3,0,0,2,0,0,0,0,0,0,1\}\\
248402&$np$&6619&28269.66&\{6498,68,17,8,5,4,0,3,0,0,2,0,0,0,0,0,0,1\}\\
253562&$np$&6732&28811.71&\{6610,69,17,8,5,4,0,3,0,0,2,0,0,0,0,0,0,1\}\\
254010&$np$&6744&28866.88&\{6621,70,16,9,5,4,0,3,0,0,2,0,0,0,0,0,0,1\}\\
278848&$np$&7250&31306.82&\{7124,72,16,9,5,4,0,3,0,0,2,0,0,0,0,0,0,0,1\}\\
295125&$n/q$&7528&32660.65&\{7401,75,16,8,5,4,3,0,0,0,0,2,0,0,0,0,0,1\}\\
295539&$np$&7535&32688.82&\{7407,74,16,9,5,4,3,0,0,0,2,0,0,0,0,0,0,0,0,1\}\\
296975&$np$&7557&32787.94&\{7427,74,17,8,5,4,0,3,0,0,0,2,0,0,0,0,0,0,0,1\}\\
300923&$np$&7615&33076.79&\{7486,75,16,9,5,4,3,0,0,0,2,0,0,0,0,0,0,0,0,1\}\\
319878&$n/q$&7909&34515.15&\{7779,77,17,8,5,4,3,0,0,0,0,2,0,0,0,0,0,1\}\\
326307&$np$&8002&34953.32&\{7868,76,18,9,5,4,0,3,0,0,0,2,0,0,0,0,0,0,0,1\}\\
343266&$n/q$&8250&36180.2&\{8117,78,18,9,5,4,3,0,0,0,2,0,0,0,0,0,0,0,1\}\\
352951&$np$&8397&36894.53&\{8262,78,17,9,6,4,0,3,0,0,2,0,0,0,0,0,0,0,0,1\}\\
369049&$n/q$&8663&38202.89&\{8526,80,17,9,6,4,0,3,0,0,2,0,0,0,0,0,0,0,0,1\}\\
383293&$np$&8934&39542.44&\{8796,80,19,9,5,4,0,3,0,0,0,2,0,0,0,0,0,0,1\}\\
399761&$n/q$&9301&41368.57&\{9163,82,18,9,5,4,3,0,0,0,0,2,0,0,0,0,0,0,1\}\\
402302&$np$&9363&41659.79&\{9221,82,18,9,6,4,0,3,0,0,0,2,0,0,0,0,0,0,0,0,1\}\\
410574&$np$&9552&42603.62&\{9410,82,19,9,6,4,0,3,0,0,0,2,0,0,0,0,0,0,0,1\}\\
425173&$np$&9873&44202.21&\{9729,84,18,9,6,4,0,3,0,0,0,0,2,0,0,0,0,0,0,1\}\\
429279&$n/q$&9951&44605.78&\{9810,84,18,8,6,4,0,3,0,0,0,2,0,0,0,0,0,0,1\}\\
437688&$n/q$&10120&45450.94&\{9978,85,18,8,6,4,0,3,0,0,0,2,0,0,0,0,0,0,1\}\\
460526&$np$&10630&47995.87&\{10482,87,19,9,6,4,0,3,0,0,0,0,2,0,0,0,0,0,0,1\}\\
466844&$np$&10765&48674.7&\{10616,88,19,9,6,4,0,3,0,0,0,0,2,0,0,0,0,0,0,1\}\\
469711&$np$&10828&48994.26&\{10679,88,20,9,6,4,0,3,0,0,0,2,0,0,0,0,0,0,0,1\}\\
475900&$np$&10961&49664.16&\{10811,89,20,9,6,4,0,3,0,0,0,2,0,0,0,0,0,0,0,1\}\\
499227&$np$&11473&52248.58&\{11320,91,19,9,6,4,0,3,0,0,0,0,2,0,0,0,0,0,0,0,1\}\\
507295&$np$&11652&53159.09&\{11499,91,19,9,6,0,4,3,0,0,0,0,2,0,0,0,0,0,0,1\}\\
509905&$np$&11705&53429.78&\{11552,91,20,9,6,4,0,3,0,0,0,0,2,0,0,0,0,0,0,1\}\\
515952&$np$&11838&54103.62&\{11684,92,19,9,6,0,4,3,0,0,0,0,2,0,0,0,0,0,0,1\}\\
522227&$np$&11985&54854.01&\{11831,92,20,9,6,4,0,3,0,0,0,0,2,0,0,0,0,0,0,1\}\\
523311&$np$&12011&54982.68&\{11856,93,19,9,6,4,0,3,0,0,0,0,2,0,0,0,0,0,0,0,1\}\\
535241&$np$&12290&56405.62&\{12134,94,19,9,6,0,4,3,0,0,0,2,0,0,0,0,0,0,0,0,1\}\\
546956&$np$&12548&57723.99&\{12391,94,21,9,6,4,0,3,0,0,0,0,2,0,0,0,0,0,0,1\}\\
557554&$np$&12782&58920.39&\{12624,95,20,9,6,0,4,3,0,0,0,0,2,0,0,0,0,0,0,1\}\\
557563&$np$&12782&58921.4&\{12624,95,21,9,6,4,0,3,0,0,0,0,2,0,0,0,0,0,0,1\}\\
564495&$np$&12930&59677.33&\{12771,96,20,9,6,4,0,3,0,0,0,0,2,0,0,0,0,0,0,0,1\}\\
580497&$np$&13263&61389.04&\{13103,97,21,9,6,4,0,3,0,0,0,2,0,0,0,0,0,0,0,0,1\}\\
588333&$np$&13433&62261.14&\{13272,98,20,9,6,4,0,3,0,0,0,0,2,0,0,0,0,0,0,0,1\}\\
640349&$np$&14657&68593.7&\{14494,99,21,10,6,0,4,3,0,0,0,2,0,0,0,0,0,0,0,1\}\\
641730&$np$&14687&68749.71&\{14524,99,22,9,6,4,0,3,0,0,0,0,2,0,0,0,0,0,0,1\}\\
720552&$np$&16235&76816.12&\{16065,105,21,10,6,0,4,3,0,0,0,2,0,0,0,0,0,0,0,0,1\}\\
773313&$np$&17205&81904.05&\{17030,108,22,10,6,0,4,3,0,0,0,0,2,0,0,0,0,0,0,0,1\}\\
783559&$np$&17384&82842.56&\{17207,110,22,9,6,0,4,3,0,0,0,0,2,0,0,0,0,0,0,0,0,1\}\\
785828&$np$&17421&83041.6&\{17245,109,22,10,6,0,4,3,0,0,0,0,2,0,0,0,0,0,0,0,1\}\\
805294&$np$&17750&84778.47&\{17573,110,23,9,6,0,4,3,0,0,0,0,2,0,0,0,0,0,0,0,1\}\\
861970&$n/q$&18828&90497.56&\{18649,114,23,10,6,0,4,3,0,0,0,2,0,0,0,0,0,0,1\}\\
869007&$np$&18969&91236.61&\{18788,114,22,10,6,0,4,0,3,0,0,0,2,0,0,0,0,0,0,1\}\\
877665&$np$&19125&92059.61&\{18942,115,23,9,6,0,4,3,0,0,0,0,2,0,0,0,0,0,0,0,0,1\}\\
881619&$n/q$&19192&92432.54&\{19012,115,23,10,6,0,4,3,0,0,0,2,0,0,0,0,0,0,1\}\\
886325&$np$&19283&92898.85&\{19099,116,23,9,6,0,4,3,0,0,0,0,2,0,0,0,0,0,0,0,0,1\}\\
902266&$np$&19552&94330.92&\{19367,117,23,9,6,0,4,3,0,0,0,0,2,0,0,0,0,0,0,0,0,1\}\\
922585&$np$&19894&96151.02&\{19707,118,23,10,6,0,4,3,0,0,0,0,2,0,0,0,0,0,0,0,0,1\}\\
936235&$np$&20122&97367.69&\{19934,119,23,10,6,0,4,3,0,0,0,0,2,0,0,0,0,0,0,0,0,1\}\\
952223&$np$&20377&98736.94&\{20190,118,24,10,6,0,4,3,0,0,0,0,2,0,0,0,0,0,0,0,1\}\\
\end{longtable}
\bibliography{main}{}

\begin{thebibliography}{1}

\bibitem{alaoglu_erdos}
L.~Alaoglu and P.~Erd\H{o}s.
\newblock On highly composite and similar numbers.
\newblock {\em Transactions of the American Mathematical Society},
  56(3):448--469, 1944.

\bibitem{noe}
T.~D. Noe.
\newblock First 1000000 superabundant numbers, 2009.
\newblock Retrieved from http://www.sspectra.com/math/SAN\_1000000.zip.

\end{thebibliography}
\bibliographystyle{plain}

\end{document}